\documentclass[12pt]{amsart}
\usepackage[top=3.6cm,bottom=3.6cm,left=3.8cm,right=3.8cm]{geometry}
\usepackage{amsmath,amssymb,mathtools}
\usepackage{booktabs}

\newcommand{\C}{\mathbb C}
\newcommand{\Q}{\mathbb Q}
\newcommand{\R}{\mathbb R}
\newcommand{\Z}{\mathbb Z}
\newcommand{\F}{\mathbb F}
\newcommand{\Hh}{\mathbb H}
\newcommand{\SL}{\mathrm{SL}_2(\Z)}

\newcommand{\dd}{\,d}

\newtheorem{theorem}{Theorem}[section]
\newtheorem{proposition}[theorem]{Proposition}
\newtheorem{corollary}[theorem]{Corollary}
\newtheorem{lemma}[theorem]{Lemma}

\numberwithin{equation}{section}

\begin{document}

\title[Jacobi pole divisors and supersingular lifts]{Reducible modular differential equations, Jacobi pole divisors, and supersingular lifts}

\author{Khalil Besrour}
\author{Hicham Saber}
\author{Abdellah Sebbar}

\address{Department of Mathematics and Statistics, University of Ottawa, Ottawa, Ontario K1N 6N5, Canada}
\email{kbesr067@uottawa.ca}
\email{asebbar@uottawa.ca}
\address{Department of Mathematics, College of Science, University of Ha'il, Ha'il 55473, Kingdom of Saudi Arabia}
\email{hi.saber@uoh.edu.sa}

\subjclass[2020]{11F03, 11F11, 33C45, 34M05}
\keywords{Modular differential equations, Schwarzian derivative, Jacobi polynomials, Kaneko-Zagier forms, supersingular invariants, singular moduli}

\begin{abstract}
We study the reducible parameter values of a second-order modular differential equation for the full modular group and show that their pole data is controlled by a single polynomial linking three a priori different structures. The residue conditions for the associated weight two form are exactly the Stieltjes equilibrium equations for four shifted Jacobi families. This determines the poles uniquely and produces a canonical divisor on $X(1)$ supported on the elliptic arc. After the standard eta-normalization, with $\ell=6r$, the differential equation takes the Kaneko-Zagier form, and the same Jacobi polynomial is precisely the polynomial factor in its distinguished modular solution. When $\ell$ is prime, this polynomial is an $\ell$-integral characteristic-zero lift of the non-elliptic supersingular polynomial. Consequently, its splitting field over $\Q_\ell$ is unramified of degree at most two, with Frobenius cycle structure determined by the supersingular $j$-invariants. Complete splitting occurs exactly when every supersingular $j$-invariant in characteristic $\ell$ is defined over $\F_\ell$, equivalently when the Fricke quotient $X_0^+(\ell)$ has genus zero. These are Ogg's primes, which are exactly the prime divisors of the order of the Monster. A common singular Sturm-Liouville equation further gives strict interlacing for consecutive reducible parameters. Finally, we prove that the canonical representatives of the poles are transcendental although their $j$-invariants are algebraic, and determine their limiting distribution and endpoint scales on the modular arc.
\end{abstract}

\maketitle

\section{Introduction}

Consider the modular differential equation
\begin{equation}\label{eq:mde-intro}
y''(z)+\pi^2r^2E_4(z)y(z)=0,
\end{equation}
where $E_4$ is the normalized Eisenstein series of weight $4$ for $\SL$ and $r\in\Q$. If $y_1,y_2$ are independent solutions, then their quotient satisfies the Schwarzian equation
\begin{equation}\label{eq:sch-intro}
\{h,z\}=2\pi^2r^2E_4(z).
\end{equation}
The modular transformation law of $E_4$ induces a two-dimensional representation on the solution space of \eqref{eq:mde-intro}, while the quotient of two independent solutions is equivariant for the corresponding projective action \cite{forum,ram1,jmaa,ram2}. The reducible case is especially rigid: the associated representation is reducible if and only if
\[
r=\frac{m}{6},\qquad \gcd(m,6)=1.
\]

For these parameters, the construction of \cite{ram1,revista} expresses $h'$ as a meromorphic modular form of weight $2$ with a character. Its non-elliptic poles are encoded by values $x_i$ of the normalized modular invariant, and the vanishing of the residues is equivalent to a finite algebraic system
\begin{equation}\label{eq:system-intro}
\frac{a}{x_i}+\frac{b}{x_i-1}+\sum_{j\ne i}\frac{12}{x_i-x_j}=0,
\qquad 1\le i\le k,
\end{equation}
for one of four pairs $(a,b)$. These systems were originally used as an existence mechanism for the modular differential equation. The point of the present paper is that they carry substantially more structure. They determine the pole divisor uniquely and connect the reducible Schwarzian problem simultaneously with Jacobi orthogonal polynomials, the Kaneko-Zagier equation, and supersingular arithmetic.

After normalization, \eqref{eq:system-intro} is exactly the Stieltjes equilibrium system for a shifted Jacobi polynomial with
\[
\alpha\in\left\{-\frac13,\frac13\right\},
\qquad
\beta\in\left\{-\frac12,\frac12\right\}.
\]
The parameters $1/3$ and $1/2$ reflect the stabilizer orders of the two elliptic points of $X(1)$. It follows that the residue equations have a unique unordered solution, that all pole invariants lie in $(0,1)$, and that the corresponding points of the modular curve have canonical representatives on the elliptic arc from $\rho=e^{2\pi i/3}$ to $i$. Thus a residue-cancellation problem arising from a differential equation produces a distinguished Jacobi divisor on $X(1)$, together with its electrostatic interpretation and exact symmetric identities.

The same polynomial has a second, arithmetic interpretation. Put $\ell=6r$. After the standard eta-normalization, equation \eqref{eq:mde-intro} takes the Kaneko-Zagier form. The Jacobi polynomial selected by the residue conditions then agrees, up to the explicit normalization dictated by the change from $J$ to $j$, with the polynomial factor in the distinguished Kaneko-Zagier modular solution. Specialized Jacobi polynomials already occur in the Kaneko-Zagier theory \cite{kz}, and Baba and Granath placed orthogonal systems of modular forms and supersingular reduction in a related framework \cite{babagranath}. The point here is an exact identification: the polynomial governing the poles in the reducible Schwarzian construction is the same polynomial governing the zero divisor in the Kaneko-Zagier solution. Thus the analytic pole problem and the arithmetic modular-form construction meet at a single canonical divisor on $X(1)$.

Write
\[
\ell=6r=12k+e,
\qquad e\in\{1,5,7,11\}.
\]
When $\ell$ is prime, the Kaneko-Zagier congruence gives an arithmetic interpretation of the pole polynomial as an $\ell$-integral characteristic-zero lift of the non-elliptic supersingular polynomial \cite{kz,cullinanhajir}. This allows local information about the pole field to be read directly from the supersingular locus. We show that its splitting field over $\Q_\ell$ is unramified and has degree at most two, and that the Frobenius cycle type is determined by the supersingular invariants defined over $\F_\ell$. In particular, complete splitting over $\Q_\ell$ occurs exactly when every supersingular $j$-invariant in characteristic $\ell$ is defined over $\F_\ell$. By Ogg's theorem, this is equivalent to the genus-zero condition
\[
\operatorname{genus}X_0^+(\ell)=0,
\]
for the Fricke quotient $X_0^+(\ell)=X_0(\ell)/\langle w_\ell\rangle$. The primes satisfying this condition are Ogg's primes; they are exactly the prime divisors of the order of the Monster \cite{ogg,duncano}. Thus the appearance of moonshine in the present problem is concrete and limited: both constructions detect the same supersingular divisor, and Ogg's criterion translates its field of definition into a complete splitting criterion for the pole polynomial. The irreducibility and Galois results of Cullinan and Hajir likewise transfer to the same pole fields \cite{cullinanhajir}.

The Jacobi description also reveals a global ordering that is not visible in the original residue equations. After multiplication by the appropriate endpoint factors, all four Jacobi families satisfy one singular Sturm-Liouville differential expression, with the spectral parameter determined by $r$. This gives a uniform comparison of the pole divisors as the reducible parameter increases. In particular, the pole configurations for every pair of consecutive reducible parameters strictly interlace on the elliptic arc. The argument includes the mixed transition between the residue classes $e=5$ and $e=7$, where the two Jacobi parameters move in opposite directions, as well as the transition from $e=11$ at degree $k$ to $e=1$ at degree $k+1$. This complements the parameter-shift interlacing results for Jacobi polynomials in \cite{driver}.

Finally, the pole divisor exhibits a sharp distinction between its arithmetic and analytic realizations. Its $j$-invariants are algebraic, but the canonical representatives on the elliptic arc are transcendental. The proof uses Schneider's theorem together with the Galois structure of singular moduli \cite{schneider,cox}. As the degree tends to infinity within any one of the four families, the pole divisors converge to the pullback of the arcsine distribution on $[0,1]$ \cite{sze}. Their extreme points approach $\rho$ and $i$ at rates $k^{-2/3}$ and $k^{-1}$, respectively. The different exponents are the geometric trace of the ramification orders $3$ and $2$ of the modular invariant at the two elliptic points.

The paper develops these connections in the order suggested above. We first derive and solve the residue system in Jacobi form. We then identify the resulting polynomial with the Kaneko-Zagier factor and study its supersingular reduction and local arithmetic. The canonical pole divisor and its interlacing properties are treated next, followed by explicit solutions, the transcendence result, and the asymptotic geometry of the pole configurations.

\section{The reducible modular construction}

We write $q=e^{2\pi iz}$ and use the normalized Eisenstein series
\[
E_4(z)=1+240\sum_{n\ge1}\sigma_3(n)q^n,
\qquad
E_6(z)=1-504\sum_{n\ge1}\sigma_5(n)q^n.
\]
The discriminant and the Dedekind eta-function satisfy
\[
\Delta(z)=\eta(z)^{24}=\frac{E_4(z)^3-E_6(z)^2}{1728}.
\]
We use the normalized Hauptmodul
\begin{equation}\label{eq:J}
J(z)=\frac{j(z)}{1728}
=\frac{E_4(z)^3}{E_4(z)^3-E_6(z)^2}.
\end{equation}
Thus
\[
J(\rho)=0,\qquad J(i)=1,\qquad J(i\infty)=\infty.
\]

For a meromorphic function $h$, its Schwarzian derivative is
\[
\{h,z\}=\left(\frac{h''}{h'}\right)'
-\frac12\left(\frac{h''}{h'}\right)^2.
\]
The following standard correspondence will be used throughout.

\begin{lemma}\label{lem:schwarzian}
Let $R$ be meromorphic on a domain. If $y_1,y_2$ are independent solutions of
\[
y''+R(z)y=0,
\]
then $h=y_2/y_1$ satisfies $\{h,z\}=2R(z)$. Conversely, if $h$ is locally univalent and $\{h,z\}=2R(z)$, then
\[
y_1=\frac1{\sqrt{h'}},\qquad y_2=\frac{h}{\sqrt{h'}}
\]
form a fundamental system.
\end{lemma}

The equivariant interpretation of the Schwarzian equation is recalled from \cite{forum}. If $\{h,z\}$ is a weight $4$ modular form, then $h$ transforms projectively under $\SL$. In the reducible case, the derivative has a simpler description.

The following result is proved in \cite{forum}.

\begin{theorem}\label{thm:reducible}
A meromorphic solution of
\[
\{h,z\}=2\pi^2r^2E_4(z)
\]
has reducible associated representation if and only if $h'$ is a meromorphic modular form of weight $2$ with a character. This occurs if and only if
\[
r=\frac{m}{6},\qquad \gcd(m,6)=1.
\]
\end{theorem}

For the construction of \cite{revista}, let
\[
\varepsilon,\delta\in\{0,2\},\qquad k\ge0,
\]
and let $w_1,\dots,w_k$ be pairwise inequivalent non-elliptic points of $\Hh$. Set
\begin{equation}\label{eq:ansatz-old}
F(z)=
\frac{\eta(z)^{4+8\varepsilon+12\delta}}
{E_4(z)^{\varepsilon}E_6(z)^{\delta}}
\prod_{j=1}^k\frac1{\bigl(J(z)-J(w_j)\bigr)^2}.
\end{equation}
Then $F$ is a meromorphic form of weight $2$ with a character. Its non-elliptic poles are double poles at the $\SL$-orbits of the $w_j$. When $\varepsilon=2$ or $\delta=2$, there is also a double pole at $\rho$ or $i$, respectively.

The residue calculation of \cite{revista} gives the following. Put
\[
x_j=J(w_j).
\]
The residues at the non-elliptic poles vanish if and only if
\begin{equation}\label{eq:system-ab}
\frac{a}{x_i}+\frac{b}{x_i-1}
+12\sum_{j\ne i}\frac1{x_i-x_j}=0,
\qquad 1\le i\le k,
\end{equation}
where
\begin{equation}\label{eq:ab}
a=4+2\varepsilon,
\qquad
b=3+3\delta.
\end{equation}
The residues at the elliptic points vanish as well. Once these residue conditions hold, $F$ has a global meromorphic primitive $h$ on $\Hh$, and
\begin{equation}\label{eq:r}
r=2k+\frac{1+2\varepsilon+3\delta}{6}.
\end{equation}
Thus the four choices give
\[
(\varepsilon,\delta)=(0,0),(2,0),(0,2),(2,2)
\quad\Longleftrightarrow\quad
6r\equiv1,5,7,11\pmod{12},
\]
respectively.

\section{The residue system and its Jacobi solution}

Define
\begin{equation}\label{eq:alphabeta}
\alpha=\frac{\varepsilon-1}{3},
\qquad
\beta=\frac{\delta-1}{2}.
\end{equation}
Hence
\[
\alpha\in\left\{-\frac13,\frac13\right\},
\qquad
\beta\in\left\{-\frac12,\frac12\right\}.
\]
Using \eqref{eq:ab}, system \eqref{eq:system-ab} becomes
\begin{equation}\label{eq:stieltjes}
\frac{\alpha+1}{x_i}
+\frac{\beta+1}{x_i-1}
+2\sum_{j\ne i}\frac1{x_i-x_j}=0,
\qquad 1\le i\le k.
\end{equation}
This is the Stieltjes system associated with a Jacobi weight on $(0,1)$.

Let $P_k^{(A,B)}$ denote the classical Jacobi polynomial. We define the monic shifted polynomial
\begin{equation}\label{eq:Qdef}
Q_k^{\alpha,\beta}(x)
=
\frac{k!}{(k+\alpha+\beta+1)_k}
P_k^{(\beta,\alpha)}(2x-1),
\end{equation}
where $(u)_k$ is the Pochhammer symbol. The normalization in \eqref{eq:Qdef} is chosen so that the leading coefficient is $1$.

The polynomial $Q_k^{\alpha,\beta}$ satisfies
\begin{equation}\label{eq:jacobi-ode}
x(1-x)Q''+
\bigl[\alpha+1-(\alpha+\beta+2)x\bigr]Q'
+k(k+\alpha+\beta+1)Q=0.
\end{equation}
Equivalently,
\begin{equation}\label{eq:hyp}
{}_2F_1\!\left(-k,k+\alpha+\beta+1;\alpha+1;x\right)
=\frac{k!}{(\alpha+1)_k}
P_k^{(\alpha,\beta)}(1-2x).
\end{equation}

\begin{theorem}\label{thm:unique-config}
For each $k\ge1$ and each of the four pairs $(\alpha,\beta)$ in \eqref{eq:alphabeta}, consider system \eqref{eq:stieltjes} on the configuration space of tuples $(x_1,\dots,x_k)\in(\C\setminus\{0,1\})^k$ with $x_i\ne x_j$ for $i\ne j$. It has a unique solution up to permutation. If
\[
0<x_{k,1}<\cdots<x_{k,k}<1
\]
are the zeros of $Q_k^{\alpha,\beta}$, then every solution of \eqref{eq:stieltjes} is a permutation of
\[
(x_{k,1},\dots,x_{k,k}).
\]
\end{theorem}

\begin{proof}
Since $\alpha,\beta>-1$, the Jacobi polynomial in \eqref{eq:Qdef} has $k$ simple zeros in $(0,1)$. At a zero $x_{k,i}$ we have
\[
\frac{Q_k''(x_{k,i})}{Q_k'(x_{k,i})}
=2\sum_{j\ne i}\frac1{x_{k,i}-x_{k,j}}.
\]
Evaluating \eqref{eq:jacobi-ode} at $x_{k,i}$ gives \eqref{eq:stieltjes}. Thus the Jacobi zeros provide a solution.

Conversely, suppose that $(x_1,\dots,x_k)$ is a solution of \eqref{eq:stieltjes} with distinct coordinates in $\C\setminus\{0,1\}$, and set
\[
P(x)=\prod_{j=1}^k(x-x_j).
\]
At each root $x_i$, equation \eqref{eq:stieltjes} gives
\[
x_i(1-x_i)P''(x_i)
+\bigl[\alpha+1-(\alpha+\beta+2)x_i\bigr]P'(x_i)=0.
\]
Consequently the polynomial
\[
R(x)=x(1-x)P''(x)
+\bigl[\alpha+1-(\alpha+\beta+2)x\bigr]P'(x)
\]
is divisible by $P$. Comparing leading coefficients gives
\[
R(x)=-k(k+\alpha+\beta+1)P(x),
\]
so $P$ satisfies \eqref{eq:jacobi-ode}.

There is only one monic polynomial solution of degree $k$ to \eqref{eq:jacobi-ode}. Indeed, if two such solutions existed, their difference would have degree $n<k$ and would satisfy the same equation. Comparing the highest-degree term would give
\[
n(n+\alpha+\beta+1)=k(k+\alpha+\beta+1),
\]
which is impossible because $\alpha+\beta+1>0$. Hence $P=Q_k^{\alpha,\beta}$.
\end{proof}

\subsection{Electrostatic characterization}

The Stieltjes system has the standard electrostatic interpretation. Place $k$ unit charges on $(0,1)$ with logarithmic pair potential
\[
-\log|x_i-x_j|,
\]
and fixed endpoint charges
\[
q_0=\frac{\alpha+1}{2},
\qquad
q_1=\frac{\beta+1}{2}.
\]
Both are positive because $\alpha,\beta>-1$. The external potential is
\begin{equation}\label{eq:external-potential}
V_{\alpha,\beta}(x)
=-q_0\log x-q_1\log(1-x).
\end{equation}
The residue system is therefore the equilibrium condition for the classical logarithmic electrostatic model associated with the Jacobi weight.

Let
\[
\mathcal C_k=\{(x_1,\dots,x_k)\in\R^k:0<x_1<\cdots<x_k<1\}.
\]
Its total potential energy is
\begin{equation}\label{eq:energy}
\mathcal E(x_1,\dots,x_k)
=\sum_{i=1}^k V_{\alpha,\beta}(x_i)
-\sum_{1\le i<j\le k}\log(x_j-x_i),
\end{equation}
or explicitly,
\[
\mathcal E
=-\frac{\alpha+1}{2}\sum_{i=1}^k\log x_i
-\frac{\beta+1}{2}\sum_{i=1}^k\log(1-x_i)
-\sum_{1\le i<j\le k}\log(x_j-x_i).
\]
Equation~\eqref{eq:energy} is the classical Stieltjes electrostatic model for the Jacobi weight; see \cite{steinerberger,sze,forrester}.

The force on the $i$th mobile charge is
\begin{equation}\label{eq:force}
\mathcal F_i
=-\frac{\partial\mathcal E}{\partial x_i}
=
\frac{\alpha+1}{2x_i}
-\frac{\beta+1}{2(1-x_i)}
+\sum_{j\ne i}\frac1{x_i-x_j}.
\end{equation}
The first two terms are the endpoint forces and the last is the mutual logarithmic repulsion. Hence,
\[
\mathcal F_i=0,\qquad 1\le i\le k,
\]
is exactly the residue system \eqref{eq:stieltjes}, after multiplication by $2$.

\begin{proposition}\label{prop:energy}
The zero vector $(x_{k,1},\dots,x_{k,k})$ is the unique global minimizer of $\mathcal E$ on $\mathcal C_k$. In the electrostatic model it is therefore the unique stable equilibrium configuration.
\end{proposition}

\begin{proof}
The force-balance equations \eqref{eq:force} are \eqref{eq:stieltjes}. Moreover,
\[
\begin{split}
v^T(\operatorname{Hess}\mathcal E)v
={}&\frac{\alpha+1}{2}\sum_{i=1}^k\frac{v_i^2}{x_i^2}
+\frac{\beta+1}{2}\sum_{i=1}^k\frac{v_i^2}{(1-x_i)^2}\\
&+\sum_{1\le i<j\le k}\frac{(v_i-v_j)^2}{(x_i-x_j)^2}.
\end{split}
\]
Since $\alpha,\beta>-1$, this quadratic form is positive for $v\ne0$, so $\mathcal E$ is strictly convex.

The energy diverges to $+\infty$ when two charges collide or when a mobile charge approaches an endpoint. Hence there is a unique global equilibrium in $\mathcal C_k$. Theorem~\ref{thm:unique-config} identifies it with the Jacobi zero vector.
\end{proof}

\subsection{Exact identities for the pole invariants}

The Jacobi description gives the following symmetric identities.

\begin{proposition}\label{prop:symmetric}
Let $x_{k,1},\dots,x_{k,k}$ be the roots of $Q_k^{\alpha,\beta}$. Then
\begin{align}
\sum_{j=1}^k x_{k,j}
&=\frac{k(k+\alpha)}{2k+\alpha+\beta},\label{eq:sum}\\
\sum_{j=1}^k\frac1{x_{k,j}}
&=\frac{k(k+\alpha+\beta+1)}{\alpha+1},\label{eq:recipsum0}\\
\sum_{j=1}^k\frac1{1-x_{k,j}}
&=\frac{k(k+\alpha+\beta+1)}{\beta+1},\label{eq:recipsum1}\\
\prod_{j=1}^k x_{k,j}
&=\frac{(\alpha+1)_k}{(k+\alpha+\beta+1)_k},\label{eq:prod0}\\
\prod_{j=1}^k(1-x_{k,j})
&=\frac{(\beta+1)_k}{(k+\alpha+\beta+1)_k}.\label{eq:prod1}
\end{align}
\end{proposition}

\begin{proof}
Formula \eqref{eq:sum} follows by comparing the two highest coefficients in the hypergeometric expression \eqref{eq:hyp}. The logarithmic derivative at $x=0$ gives
\[
-\sum_{j=1}^k\frac1{x_{k,j}}
=\frac{Q_k'(0)}{Q_k(0)}.
\]
Putting $x=0$ in \eqref{eq:jacobi-ode} gives
\[
(\alpha+1)Q_k'(0)+k(k+\alpha+\beta+1)Q_k(0)=0,
\]
which yields \eqref{eq:recipsum0}. Applying the same argument to $Q_k(1-x)$ gives \eqref{eq:recipsum1}. The endpoint values of Jacobi polynomials yield
\[
Q_k(0)=(-1)^k\frac{(\alpha+1)_k}{(k+\alpha+\beta+1)_k},
\qquad
Q_k(1)=\frac{(\beta+1)_k}{(k+\alpha+\beta+1)_k},
\]
which are equivalent to \eqref{eq:prod0} and \eqref{eq:prod1}.
\end{proof}

Since $\alpha,\beta\in\Q$, the polynomial $Q_k^{\alpha,\beta}$ belongs to $\Q[x]$; hence every pole invariant $x_{k,j}$ is algebraic.

\section{The Kaneko-Zagier form and supersingular reduction}\label{sec:kz}

The four Jacobi families of Section~3 are the specialized Jacobi families occurring in the Kaneko-Zagier theory. Kaneko and Zagier obtain these polynomials from modular forms satisfying a second-order differential equation and relate their reductions to supersingular polynomials \cite{kz}. Brillhart and Morton independently studied the same Jacobi lift in connection with supersingular polynomials and class-number phenomena \cite{brillhartmorton}; the coincidence with the Kaneko-Zagier lift is also discussed in \cite{cullinanhajir}. Baba and Granath subsequently extended this perspective to orthogonal systems of modular forms with supersingular reduction \cite{babagranath}. In the present setting, the polynomial first arises from the residue conditions of the reducible Schwarzian construction and is then identified with the Kaneko-Zagier polynomial.

Our normalized Hauptmodul is $J=j/1728$; in this section the classical modular invariant is therefore $j=1728J$.

Put
\begin{equation}\label{eq:ell-e}
\ell=6r=12k+e,
\qquad
e=1+2\varepsilon+3\delta,
\end{equation}
and set
\begin{equation}\label{eq:lambdamu}
\lambda=\varepsilon-1,
\qquad
\mu=\delta-1.
\end{equation}
Then $\lambda,\mu\in\{\pm1\}$,
\[
e=6+2\lambda+3\mu,
\qquad
\alpha=\frac{\lambda}{3},
\qquad
\beta=\frac{\mu}{2}.
\]
The notation used in this section is summarized by
\begin{center}
\begin{tabular}{c l}
\toprule
symbol & meaning \\
\midrule
$r$ & parameter of the normal-form modular differential equation \\
$\ell$ & numerator $6r$, so that $\ell=12k+e$ \\
$e$ & residue class $1,5,7,$ or $11$ modulo $12$ \\
$(\varepsilon,\delta)$ & endpoint data in $\{0,2\}^2$ with $e=1+2\varepsilon+3\delta$ \\
$(\alpha,\beta)$ & Jacobi parameters $((\varepsilon-1)/3,(\delta-1)/2)$ \\
$Q_k^{\alpha,\beta}$ & monic shifted Jacobi pole polynomial in $J$ \\
$\widetilde{\mathcal F}_{\ell-1}$ & Kaneko-Zagier polynomial factor in the classical $j$-variable \\
$\mathcal P_\ell$ & monic pole polynomial in the classical $j$-variable \\
\bottomrule
\end{tabular}
\end{center}
The four choices of $(\varepsilon,\delta)$ correspond to $e=1,5,7,11$.

Let
\[
D=q\frac{d}{dq}=\frac{1}{2\pi i}\frac{d}{dz},
\qquad
\vartheta_s=D-\frac{s}{12}E_2.
\]
For $\ell\ge5$ coprime to $6$, Kaneko and Zagier consider the distinguished modular form $\mathcal F_{\ell-1}$ of weight $\ell-1$ satisfying
\begin{equation}\label{eq:KZeq}
\vartheta_{\ell+1}\vartheta_{\ell-1}\mathcal F_{\ell-1}
-\frac{\ell^2-1}{144}E_4\mathcal F_{\ell-1}=0.
\end{equation}
See \cite{kz}; the same equation and its modular solutions are also studied in \cite{nakaya}. The case $\ell=1$ is the weight zero endpoint and is omitted below. This normalization is standard; see also \cite{ittersum}. The calculation fixes the normalization of the parameter.

\begin{proposition}\label{prop:eta-normalization}
If $f$ satisfies \eqref{eq:KZeq} and
\[
f=\eta^{2\ell}y,
\]
then
\begin{equation}\label{eq:normal-eta}
y''+\pi^2\left(\frac{\ell}{6}\right)^2E_4y=0.
\end{equation}
Hence the normal-form parameter is $r=\ell/6$.
\end{proposition}

\begin{proof}
The identities
\[
D\log\eta=\frac{E_2}{24},
\qquad
DE_2=\frac{E_2^2-E_4}{12}
\]
give
\[
\vartheta_{\ell+1}\vartheta_{\ell-1}(\eta^{2\ell}y)
=\eta^{2\ell}\left(D^2y-\frac{E_4}{144}y\right).
\]
Substitution in \eqref{eq:KZeq} yields
\[
D^2y-\frac{\ell^2}{144}E_4y=0.
\]
Since $D^2=-(4\pi^2)^{-1}d^2/dz^2$, this is \eqref{eq:normal-eta}.
\end{proof}

The modular form $\mathcal F_{\ell-1}$ has a factorization
\begin{equation}\label{eq:KZfactor}
\mathcal F_{\ell-1}(z)
=\Delta(z)^kE_4(z)^{\varepsilon/2}E_6(z)^{\delta/2}
\widetilde{\mathcal F}_{\ell-1}(j(z)),
\end{equation}
where $\widetilde{\mathcal F}_{\ell-1}\in\Q[X]$ has degree $k$. Kaneko and Zagier give the corresponding hypergeometric description in \cite[\S8, Theorem~5]{kz}, and the remark following their Theorem~6 identifies the resulting polynomial, up to rescaling, with a Jacobi polynomial. In the present normalization,
\begin{equation}\label{eq:KZjacobi}
\widetilde{\mathcal F}_{\ell-1}(X)
=1728^kP_k^{(\alpha,\beta)}\!\left(1-\frac{X}{864}\right).
\end{equation}
This formula is also the starting point of the arithmetic study in \cite{cullinanhajir}.

\begin{theorem}\label{thm:KZ-pole-poly}
For the parameters in \eqref{eq:ell-e}, the Kaneko-Zagier polynomial and the pole polynomial are related by
\begin{equation}\label{eq:KZ-Q}
\widetilde{\mathcal F}_{\ell-1}(1728x)
=(-1)^k1728^k
\frac{(k+\alpha+\beta+1)_k}{k!}
Q_k^{\alpha,\beta}(x).
\end{equation}
Consequently,
\begin{equation}\label{eq:KZ-y1}
\frac{\mathcal F_{\ell-1}(z)}{\eta(z)^{2\ell}}
=C_{k,\varepsilon,\delta}
\frac{E_4(z)^{\varepsilon/2}E_6(z)^{\delta/2}
Q_k^{\alpha,\beta}(J(z))}
{\eta(z)^{2+4\varepsilon+6\delta}}
\end{equation}
for a nonzero rational constant $C_{k,\varepsilon,\delta}$.
\end{theorem}

\begin{proof}
The Jacobi reflection identity gives
\[
P_k^{(\alpha,\beta)}(1-2x)
=(-1)^kP_k^{(\beta,\alpha)}(2x-1).
\]
Combining this with \eqref{eq:Qdef} and \eqref{eq:KZjacobi} proves \eqref{eq:KZ-Q}. Next, use $\Delta=\eta^{24}$ and
\[
2\ell=24k+2+4\varepsilon+6\delta
\]
in \eqref{eq:KZfactor}. Formula \eqref{eq:KZ-y1} follows.
\end{proof}

For every admissible value $\ell=12k+e$, define the monic pole polynomial in the classical $j$-variable by
\begin{equation}\label{eq:pole-j-poly}
\mathcal P_\ell(X)
=\prod_{j=1}^k\bigl(X-1728x_{k,j}\bigr)
=1728^kQ_k^{\alpha,\beta}\!\left(\frac{X}{1728}\right).
\end{equation}
Thus $\mathcal P_\ell$ is defined before any primality assumption on $\ell$.

\begin{proposition}\label{prop:pole-poly}
The polynomial $\mathcal P_\ell$ is monic of degree $k$ and belongs to $\Q[X]$. For $\ell\ge5$, it is the monic normalization of the Kaneko-Zagier polynomial:
\begin{equation}\label{eq:P-Ftilde}
\mathcal P_\ell(X)
=(-1)^k\frac{k!}{(k+\alpha+\beta+1)_k}
\widetilde{\mathcal F}_{\ell-1}(X).
\end{equation}
For every admissible $\ell$, it satisfies
\begin{equation}\label{eq:P-ode}
\begin{aligned}
X(1728-X)\mathcal P_\ell''(X)
&+\bigl[1728(\alpha+1)-(\alpha+\beta+2)X\bigr]
\mathcal P_\ell'(X)\\
&+k(k+\alpha+\beta+1)\mathcal P_\ell(X)=0.
\end{aligned}
\end{equation}
Its $k$ roots are simple and lie in $(0,1728)$.
\end{proposition}

\begin{proof}
The first equality in \eqref{eq:pole-j-poly} and the fact that $Q_k^{\alpha,\beta}$ is monic give the second equality. When $\ell\ge5$, formula \eqref{eq:P-Ftilde} follows from \eqref{eq:KZ-Q} after setting $X=1728x$. To obtain \eqref{eq:P-ode}, substitute
\[
Q_k^{\alpha,\beta}(x)=1728^{-k}\mathcal P_\ell(1728x)
\]
into \eqref{eq:jacobi-ode}. The roots of $Q_k^{\alpha,\beta}$ are simple and belong to $(0,1)$, so \eqref{eq:pole-j-poly} places the roots of $\mathcal P_\ell$ in $(0,1728)$.
\end{proof}

For $\ell\ge5$, the modular form $\mathcal F_{\ell-1}(z)$, its polynomial factor $\widetilde{\mathcal F}_{\ell-1}(X)$, and the monic pole polynomial $\mathcal P_\ell(X)$ will be kept distinct below. Formula \eqref{eq:P-Ftilde} is the precise normalization relating the last two.

\begin{corollary}\label{cor:residue-from-ode}
Let
\[
y(z)=
\frac{E_4(z)^{\varepsilon/2}E_6(z)^{\delta/2}
Q_k^{\alpha,\beta}(J(z))}
{\eta(z)^{2+4\varepsilon+6\delta}}.
\]
At every zero $a\in\Hh$ of $y$, the meromorphic differential $y(z)^{-2}\,dz$ has zero residue. In particular, the Jacobi roots satisfy all residue conditions of the reducible construction, including the conditions at the elliptic points when they occur.
\end{corollary}

\begin{proof}
By Proposition~\ref{prop:eta-normalization} and Theorem~\ref{thm:KZ-pole-poly}, $y$ satisfies \eqref{eq:normal-eta}. A zero of a nonzero solution of a second-order equation with holomorphic coefficient is simple. If $t=z-a$, write
\[
y(z)=c_1t+c_2t^2+O(t^3).
\]
Evaluating the differential equation at $a$ gives $y''(a)=0$, hence $c_2=0$. Therefore
\[
y(z)^{-2}=c_1^{-2}t^{-2}+O(1),
\]
with no $t^{-1}$ term.
\end{proof}

\subsection{Supersingular reduction}

Assume now that $\ell=12k+e$ is prime. Let $S_\ell(X)$ be the supersingular polynomial in characteristic $\ell$. With the notation above it has the factorization
\begin{equation}\label{eq:ss-factor}
S_\ell(X)=X^{\varepsilon/2}(X-1728)^{\delta/2}s_\ell(X),
\end{equation}
where $s_\ell\in\F_\ell[X]$ has degree $k$ and contains the non-elliptic supersingular invariants. Kaneko and Zagier prove that the polynomial factor attached to their distinguished modular solution has $\ell$-integral coefficients and reduces to the supersingular polynomial after the elliptic factors are removed; see \cite[Theorem~1]{kz}. With the normalization above, equivalently,
\begin{equation}\label{eq:KZ-cong}
\widetilde{\mathcal F}_{\ell-1}(X)\equiv s_\ell(X)\pmod\ell.
\end{equation}
See also \cite{cullinanhajir} for this normalization and its arithmetic consequences.

\begin{theorem}\label{thm:ss-lift}
If $\ell=12k+e$ is prime, then $\mathcal P_\ell$ has $\ell$-integral coefficients and
\begin{equation}\label{eq:pole-ss-cong}
\mathcal P_\ell(X)\equiv s_\ell(X)\pmod\ell.
\end{equation}
Thus $\mathcal P_\ell$ is a characteristic zero lift of the non-elliptic supersingular polynomial.
\end{theorem}

\begin{proof}
By Proposition~\ref{prop:pole-poly}, $\mathcal P_\ell$ is the monic normalization of $\widetilde{\mathcal F}_{\ell-1}$. The congruence \eqref{eq:KZ-cong} has monic right-hand side. Hence the leading coefficient of $\widetilde{\mathcal F}_{\ell-1}$ is an $\ell$-adic unit congruent to $1$ modulo $\ell$. Division by this unit preserves $\ell$-integrality and gives \eqref{eq:pole-ss-cong}.
\end{proof}

The differential equation for the pole polynomial also reduces modulo $\ell$.

\begin{corollary}\label{cor:ss-ode}
For prime $\ell=12k+e$, the non-elliptic supersingular polynomial $s_\ell$ satisfies
\begin{equation}\label{eq:ss-ode}
X(1728-X)s_\ell''(X)
+\bigl[1728(\alpha+1)-(\alpha+\beta+2)X\bigr]s_\ell'(X)
+k(k+\alpha+\beta+1)s_\ell(X)=0
\end{equation}
in $\F_\ell[X]$, where the rational coefficients are reduced modulo $\ell$.
\end{corollary}

\begin{proof}
The denominators of $\alpha$ and $\beta$ divide $6$, so they are units modulo $\ell\ge5$. Reduce \eqref{eq:P-ode} modulo $\ell$ and use Theorem~\ref{thm:ss-lift}.
\end{proof}

Let $\Sigma_\ell^*=\{\xi_1,\dots,\xi_k\}$ be the roots of $s_\ell$ in $\overline{\F}_\ell$. Evaluating \eqref{eq:ss-ode} at these roots gives the characteristic-$\ell$ Stieltjes relations.

\begin{corollary}\label{cor:ss-stieltjes}
For every $1\le i\le k$,
\begin{equation}\label{eq:ss-stieltjes}
\frac{\alpha+1}{\xi_i}
+\frac{\beta+1}{\xi_i-1728}
+2\sum_{j\ne i}\frac1{\xi_i-\xi_j}=0
\end{equation}
in $\overline{\F}_\ell$.
\end{corollary}

\begin{proof}
The polynomial $s_\ell$ has distinct roots and none of them is $0$ or $1728$. At a root $\xi_i$,
\[
\frac{s_\ell''(\xi_i)}{s_\ell'(\xi_i)}
=2\sum_{j\ne i}\frac1{\xi_i-\xi_j}.
\]
Evaluating \eqref{eq:ss-ode} at $X=\xi_i$ and dividing by $\xi_i(1728-\xi_i)s_\ell'(\xi_i)$ gives \eqref{eq:ss-stieltjes} after partial fractions.
\end{proof}

\begin{corollary}\label{cor:unramified}
For prime $\ell=12k+e$, reduction modulo $\ell$ induces a bijection between the roots of $\mathcal P_\ell$ in $\overline{\Q}_\ell$ and the non-elliptic supersingular $j$-invariants in $\overline{\F}_\ell$. The splitting field of $\mathcal P_\ell$ over $\Q_\ell$ is unramified. Every irreducible factor of $\mathcal P_\ell$ over $\Q_\ell$ has degree $1$ or $2$; equivalently, the local splitting field has residue degree at most $2$. Consequently, the local splitting field is either $\Q_\ell$ or the unique unramified quadratic extension of $\Q_\ell$.
\end{corollary}

\begin{proof}
Fix an algebraic closure $\overline{\Q}_\ell$, extend the $\ell$-adic valuation to it, and identify the residue field of its valuation ring with $\overline{\F}_\ell$. Factor
\[
s_\ell(X)=f_1(X)\cdots f_t(X)
\]
into distinct monic irreducible factors over $\F_\ell$. Since $s_\ell$ is squarefree, Hensel's lemma lifts this factorization uniquely to
\[
\mathcal P_\ell(X)=F_1(X)\cdots F_t(X)
\]
with pairwise coprime monic $F_i\in\Z_\ell[X]$ satisfying $F_i\bmod\ell=f_i$ and $\deg F_i=\deg f_i=:d_i$.

Each $F_i$ is irreducible over $\Q_\ell$. Indeed, a nontrivial factorization over $\Q_\ell$ may be made integral and monic by Gauss's lemma, and reduction would give a nontrivial factorization of the irreducible polynomial $f_i$. Let $\theta_i$ be a root of $F_i$. Then
\[
[\Q_\ell(\theta_i):\Q_\ell]=d_i.
\]
The element $\theta_i$ is integral, and its reduction has minimal polynomial $f_i$, so the residue field of $\Q_\ell(\theta_i)$ contains $\F_{\ell^{d_i}}$. Its residue degree is therefore at least $d_i$, while it cannot exceed the degree of the extension. Hence the residue degree is $d_i$ and the ramification index is $1$. Thus $\Q_\ell(\theta_i)/\Q_\ell$ is the unramified extension of degree $d_i$. All roots of $F_i$ are its Frobenius conjugates, and the compositum of these unramified extensions is unramified. Therefore the splitting field of $\mathcal P_\ell$ is unramified.

Because the reduction is squarefree, distinct roots of $\mathcal P_\ell$ have distinct reductions, and Hensel lifting gives one root above each root of $s_\ell$. Thus reduction induces the asserted bijection. Every supersingular $j$-invariant in characteristic $\ell$ belongs to $\F_{\ell^2}$; see \cite{silverman}. Therefore each $d_i$ is $1$ or $2$. Hence every irreducible factor of $\mathcal P_\ell$ over $\Q_\ell$ has degree $1$ or $2$, and the unramified splitting field has residue degree at most $2$.
\end{proof}

Let
\[
a_\ell=\#\bigl(\Sigma_\ell^*\cap\F_\ell\bigr).
\]
The integer $a_\ell$ is therefore the number of non-elliptic supersingular invariants in $\F_\ell$.

Let $h_\ell$ denote the class number of the imaginary quadratic field $\Q(\sqrt{-\ell})$. Deuring's formula for the number of supersingular invariants in $\F_\ell$ gives
\begin{equation}\label{eq:deuring-count}
a_\ell+\frac{\varepsilon+\delta}{2}
=
\begin{cases}
\frac12 h_\ell,& \ell\equiv1\pmod4,\\[2mm]
2h_\ell,& \ell\equiv3\pmod8,\\[2mm]
h_\ell,& \ell\equiv7\pmod8.
\end{cases}
\end{equation}
See \cite{morton}.

\begin{proposition}\label{prop:frob}
Let $K_\ell$ be the splitting field of $\mathcal P_\ell$ over $\Q_\ell$. Arithmetic Frobenius acts on the $k$ roots of $\mathcal P_\ell$ with cycle type
\[
1^{a_\ell}2^{(k-a_\ell)/2}.
\]
In particular,
\[
\mathcal P_\ell \text{ splits completely over }\Q_\ell
\quad\Longleftrightarrow\quad
\Sigma_\ell^*\subset\F_\ell.
\]
\end{proposition}

\begin{proof}
By Corollary~\ref{cor:unramified}, reduction identifies the roots of $\mathcal P_\ell$ with the roots of $s_\ell$, and the local splitting field is unramified. Frobenius on the lifted roots therefore reduces to
\[
\xi\longmapsto \xi^\ell
\]
on $\Sigma_\ell^*$. Every supersingular invariant lies in $\F_{\ell^2}$, so the Frobenius orbits have length one or two. The length-one orbits are precisely the $a_\ell$ elements of $\Sigma_\ell^*\cap\F_\ell$. This gives the asserted cycle type and the splitting criterion.
\end{proof}

Ogg determined the primes for which every supersingular $j$-invariant is defined over the prime field. His list is also the set of prime divisors of the order of the Monster group; see \cite{ogg,duncano}.

\begin{corollary}\label{cor:ogg}
Let $\ell\ge5$ be prime. Then
\[
\mathcal P_\ell \text{ splits completely over }\Q_\ell
\]
if and only if
\[
\ell\in
\{5,7,11,13,17,19,23,29,31,41,47,59,71\}.
\]
Equivalently,
\[
\mathcal P_\ell \text{ splits completely over }\Q_\ell
\quad\Longleftrightarrow\quad
\ell\mid |\mathbb{M}|.
\]
\end{corollary}

\begin{proof}
The omitted supersingular invariants $0$ and $1728$, when present, already lie in $\F_\ell$. Proposition~\ref{prop:frob} therefore shows that $\mathcal P_\ell$ splits over $\Q_\ell$ if and only if all supersingular $j$-invariants in characteristic $\ell$ lie in $\F_\ell$. Ogg's theorem identifies the primes with this property as
\[
2,3,5,7,11,13,17,19,23,29,31,41,47,59,71.
\]
Since $\ell\ge5$, the displayed list follows. These are exactly the prime divisors of the order of the Monster.
\end{proof}

\subsection{The Ogg primes and moonshine}

The occurrence of the prime divisors of the Monster in Corollary~\ref{cor:ogg} comes through the supersingular reduction, rather than from a direct action of the Monster on the differential equation. Ogg's theorem identifies the same set of primes through the Fricke quotient
\[
X_0^+(\ell)=X_0(\ell)/\langle w_\ell\rangle.
\]
In characteristic $\ell$, the special fiber of $X_0(\ell)$ is described by two copies of $X(1)$ meeting at the supersingular points. The Fricke involution exchanges the two components and acts as Frobenius on their intersection. This picture explains geometrically why the field of definition of the supersingular points enters the genus calculation. The precise result used here is Ogg's calculation
\begin{equation}\label{eq:ogg-geometric}
\operatorname{genus}X_0^+(\ell)=0
\quad\Longleftrightarrow\quad
\begin{gathered}
\text{all supersingular $j$-invariants in characteristic $\ell$}\\
\text{lie in $\F_\ell$.}
\end{gathered}
\end{equation}
See \cite{ogg,duncano}.

Monstrous moonshine supplies the genus-zero side of \eqref{eq:ogg-geometric}. For each prime $\ell$ dividing the order of the Monster, there is a prime-order conjugacy class whose McKay--Thompson series is a Hauptmodul for the Fricke group $\Gamma_0(\ell)^+$. Hence $X_0^+(\ell)$ has genus zero. Combined with \eqref{eq:ogg-geometric}, this gives the implication
\[
\ell\mid |\mathbb M|
\quad\Longrightarrow\quad
\Sigma_\ell^*\subset\F_\ell,
\]
and Corollary~\ref{cor:ogg} translates it into complete splitting of the pole polynomial over $\Q_\ell$.

Duncan and Ono also give a comparison at the level of the supersingular divisor. For a prime-order Monster element with moonshine group $\Gamma_0(\ell)^+$, the reduction modulo $\ell$ of a suitable $U(\ell)$-transform of its McKay--Thompson series is a linear combination of reciprocals of polynomials
\[
 j(z)-\xi,
\qquad \xi\in\Sigma_\ell^*.
\]
By Theorem~\ref{thm:ss-lift}, these are exactly the linear factors of the reduction of $\mathcal P_\ell$. The Kaneko-Zagier pole polynomial and the prime-order moonshine functions therefore detect the same non-elliptic supersingular divisor after reduction modulo $\ell$, although they arise from different characteristic-zero constructions.

No direct Monster action on the solution space of \eqref{eq:mde-intro} follows from this comparison, nor does the pole configuration alone produce an element of order $\ell$ in the Monster. The connection used here is through the common supersingular divisor and Ogg's genus-zero criterion.

Equation~\eqref{eq:deuring-count} expresses the Frobenius cycle type in Proposition~\ref{prop:frob} in terms of an imaginary quadratic class number.

Proposition~\ref{prop:symmetric} gives the following reductions.

\begin{corollary}\label{cor:ss-symmetric}
For prime $\ell=12k+e$, the following identities hold in $\overline{\F}_\ell$:
\begin{align}
\sum_{\xi\in\Sigma_\ell^*}\xi
&=1728\frac{k(k+\alpha)}{2k+\alpha+\beta},\label{eq:ss-sum}\\
\sum_{\xi\in\Sigma_\ell^*}\frac1\xi
&=\frac1{1728}\frac{k(k+\alpha+\beta+1)}{\alpha+1},\label{eq:ss-recipsum0}\\
\sum_{\xi\in\Sigma_\ell^*}\frac1{1728-\xi}
&=\frac1{1728}\frac{k(k+\alpha+\beta+1)}{\beta+1},\label{eq:ss-recipsum1}\\
\prod_{\xi\in\Sigma_\ell^*}\xi
&=1728^k\frac{(\alpha+1)_k}{(k+\alpha+\beta+1)_k},\label{eq:ss-prod0}\\
\prod_{\xi\in\Sigma_\ell^*}(1728-\xi)
&=1728^k\frac{(\beta+1)_k}{(k+\alpha+\beta+1)_k}.\label{eq:ss-prod1}
\end{align}
The rational quantities on the right are interpreted through their reductions modulo $\ell$.
\end{corollary}

\begin{proof}
By Theorem~\ref{thm:ss-lift}, the roots of $s_\ell$ are the reductions of the numbers $1728x_{k,j}$ through their monic defining polynomial. Reducing \eqref{eq:sum}, \eqref{eq:recipsum0}, \eqref{eq:recipsum1}, \eqref{eq:prod0}, and \eqref{eq:prod1} gives the formulas. The displayed denominators are units modulo $\ell$. Indeed,
\[
2k+\alpha+\beta=\frac{\ell-6}{6},
\]
while $\alpha+1$ and $\beta+1$ have denominators dividing $6$ and nonzero numerators modulo $\ell$. If $0\le t\le k-1$, then
\[
6(k+\alpha+\beta+1+t)=\ell-6(k-t),
\]
which is nonzero modulo $\ell$. Thus every factor in $(k+\alpha+\beta+1)_k$ is an $\ell$-adic unit.
\end{proof}

\subsection{Pole fields and Galois groups}

The affine change $X=1728x$ and multiplication by a nonzero rational constant do not alter the splitting field, so the results of \cite{cullinanhajir} apply directly to the pole polynomials.

\begin{proposition}\label{prop:galois-transfer}
The polynomials $Q_k^{\alpha,\beta}$ and $\widetilde{\mathcal F}_{\ell-1}$ have isomorphic splitting fields over $\Q$ and the same irreducibility behavior. In particular:
\begin{enumerate}
\item if $\widetilde{\mathcal F}_{\ell-1}$ is irreducible and
\[
\ell=p+6q
\]
with $p,q$ prime and
\[
\frac{k}{2}<q<k-2,
\]
then
\[
\operatorname{Gal}\bigl(Q_k^{\alpha,\beta}/\Q\bigr)\cong S_k;
\]
\item the following four infinite subfamilies are irreducible over $\Q$:
\begin{align*}
(\varepsilon,\delta)=(0,0),&\qquad k=2\cdot4^{\nu-1}, &&\nu\ge1,\\
(\varepsilon,\delta)=(2,0),&\qquad k=4^{\nu-1}, &&\nu\ge1,\\
(\varepsilon,\delta)=(0,2),&\qquad k=2\cdot4^{\nu-1}-1, &&\nu\ge1,\\
(\varepsilon,\delta)=(2,2),&\qquad k=4^{\nu-1}-1, &&\nu\ge2.
\end{align*}
\end{enumerate}
In every irreducible case, the $k$ non-elliptic pole invariants form a single Galois orbit.
\end{proposition}

\begin{proof}
The first assertion follows immediately from \eqref{eq:KZ-Q}. The Galois criterion and the four irreducible families are Theorems~1.6 and~1.9 of Cullinan and Hajir \cite{cullinanhajir}, translated by
\[
\ell=12k+1+2\varepsilon+3\delta.
\]
\end{proof}

\section{Canonical pole divisors on the modular curve}

Let
\begin{equation}\label{eq:arc}
\mathcal A=\{e^{i\theta}:\pi/2\le\theta\le2\pi/3\}
\end{equation}
be the elliptic arc in the standard fundamental domain, with endpoints $i$ and $\rho=e^{2\pi i/3}$. The classical modular invariant $j$ maps this arc bijectively onto $[0,1728]$; see, for example, \cite{apostol}. With the normalization \eqref{eq:J}, $J$ maps $\mathcal A$ bijectively onto $[0,1]$, with
\[
J(\rho)=0,\qquad J(i)=1.
\]
We orient $\mathcal A$ from $\rho$ to $i$.

For each zero $x_{k,j}$, let $w_{k,j}$ be the unique point in the interior of $\mathcal A$ satisfying
\begin{equation}\label{eq:wdef}
J(w_{k,j})=x_{k,j}.
\end{equation}

\begin{theorem}\label{thm:pole-divisor}
Fix $(\varepsilon,\delta)\in\{0,2\}^2$ and $k\ge1$. The non-elliptic pole divisor required by the reducible modular construction is unique on $X(1)$. Its representatives in the standard fundamental domain are 
\[
w_{k,1},\dots,w_{k,k}\in\mathcal A^{\circ}.
\]
In particular, every non-elliptic pole lies on the $\SL$-orbit of the elliptic arc.
\end{theorem}

\begin{proof}
Theorem~\ref{thm:unique-config} gives a unique unordered set of values
\[
\{x_{k,1},\dots,x_{k,k}\}\subset(0,1).
\]
Since $J$ is a Hauptmodul for $X(1)$, each value determines one point of $X(1)$. The restriction of $J$ to $\mathcal A$ gives the stated representatives.
\end{proof}

The four Jacobi families satisfy a common Sturm-Liouville equation.

For $(\varepsilon,\delta)\in\{0,2\}^2$, put
\begin{equation}\label{eq:Udef}
U_{k,\varepsilon,\delta}(x)
=
x^{\varepsilon/6}(1-x)^{\delta/4}
Q_k^{\alpha,\beta}(x),
\end{equation}
where $\alpha,\beta$ are given by \eqref{eq:alphabeta}. Set
\[
p(x)=x^{2/3}(1-x)^{1/2},
\qquad
w(x)=x^{-1/3}(1-x)^{-1/2}.
\]

\begin{proposition}\label{prop:common-sturm}
The function $U_{k,\varepsilon,\delta}$ satisfies
\begin{equation}\label{eq:common-sturm}
\bigl(pU'\bigr)'
+
\left(\frac{r^2}{4}-\frac1{144}\right)wU=0,
\end{equation}
where $r$ is given by \eqref{eq:r}. Its zeros in $(0,1)$ are the roots of $Q_k^{\alpha,\beta}$. At $x=0$ it has the behavior
\[
U_{k,\varepsilon,\delta}(x)\asymp x^{\varepsilon/6},
\]
and at $x=1$,
\[
U_{k,\varepsilon,\delta}(x)\asymp (1-x)^{\delta/4}.
\]
\end{proposition}

\begin{proof}
Substitute
\[
Q(x)=x^{-\varepsilon/6}(1-x)^{-\delta/4}U(x)
\]
into \eqref{eq:jacobi-ode}, and use \eqref{eq:alphabeta} and \eqref{eq:r}. After cancellation one obtains
\[
x(1-x)U''
+\left(\frac23-\frac76x\right)U'
+\left(\frac{r^2}{4}-\frac1{144}\right)U=0.
\]
Since
\[
\frac{p'}{w}=\frac23-\frac76x,
\qquad
\frac{p}{w}=x(1-x),
\]
this is \eqref{eq:common-sturm}. The endpoint assertions follow from the nonzero endpoint values of the relevant Jacobi polynomials.
\end{proof}

\begin{theorem}\label{thm:interlace}
For fixed $(\varepsilon,\delta)$, the pole configurations of degrees $k$ and $k+1$ strictly interlace on $\mathcal A$. More precisely,
\[
0<x_{k+1,1}<x_{k,1}<x_{k+1,2}<\cdots
<x_{k,k}<x_{k+1,k+1}<1.
\]
Equivalently, along the orientation from $\rho$ to $i$,
\[
\rho\prec w_{k+1,1}\prec w_{k,1}\prec w_{k+1,2}
\prec\cdots\prec w_{k,k}\prec w_{k+1,k+1}\prec i.
\]
The corresponding parameters in the modular differential equation differ by $2$:
\[
r_{k+1}=r_k+2.
\]
\end{theorem}

\begin{proof}
The Jacobi weight
\[
x^\alpha(1-x)^\beta\dd x
\]
is positive on $(0,1)$. The zeros of consecutive orthogonal polynomials therefore strictly interlace; see \cite{sze}. The modular statement follows from the monotonicity of $J$ on the oriented arc. Formula \eqref{eq:r} gives the last assertion.
\end{proof}

Write $x_{k,j}^{(e)}$ for the roots in the residue class $e\in\{1,5,7,11\}$, in increasing order. The classes $e=1,5,7,11$ correspond respectively to
\[
(\varepsilon,\delta)=(0,0),(2,0),(0,2),(2,2).
\]

Let
\[
W[U,V](x)=p(x)\bigl(U(x)V'(x)-V(x)U'(x)\bigr).
\]
If $U$ and $V$ satisfy \eqref{eq:common-sturm} with spectral parameters $\Lambda_U$ and $\Lambda_V$, then
\begin{equation}\label{eq:wronskian-derivative}
W[U,V]'=(\Lambda_U-\Lambda_V)wUV.
\end{equation}
Thus $W[U,V]$ is strictly decreasing on every interval on which $U,V>0$ and $\Lambda_V>\Lambda_U$.

\begin{lemma}\label{lem:endpoint-W}
Suppose the leading constants below are positive. Then:
\begin{enumerate}
\item if $U(x)=A x^{1/3}(1+O(x))$ and $V(x)=B(1+O(x))$ as $x\to0^+$, then
\[
\lim_{x\to0^+}W[U,V](x)=-\frac{AB}{3};
\]
\item if $U(x)=A(1-x)^{1/2}(1+O(1-x))$ and $V(x)=B(1+O(1-x))$ as $x\to1^-$, then
\[
\lim_{x\to1^-}W[U,V](x)=\frac{AB}{2};
\]
\item if both $U$ and $V$ have nonzero finite limits at $1$, then $W[U,V](x)\to0$ as $x\to1^-$;
\item if
\[
U(x)=A(1-x)^{1/2}(1+O(1-x)),\qquad
V(x)=B(1-x)^{1/2}(1+O(1-x)),
\]
then $W[U,V](x)\to0$ as $x\to1^-$.
\end{enumerate}
\end{lemma}

\begin{proof}
Differentiate the displayed endpoint expansions and use
\[
p(x)\sim x^{2/3}\quad(x\to0^+),
\qquad
p(x)\sim(1-x)^{1/2}\quad(x\to1^-).
\]
In the last case the leading terms in $UV'-VU'$ cancel; the remaining term is bounded after the common factor $(1-x)^{1/2}$ is removed, and multiplication by $p(x)$ gives the limit $0$.
\end{proof}

Driver, Jordaan, and Mbuyi \cite{driver} give several interlacing criteria for Jacobi polynomials of equal or adjacent degrees under continuous parameter shifts. With our convention $Q_k^{\alpha,\beta}(x)\propto P_k^{(\beta,\alpha)}(2x-1)$, their results recover the transitions $e=1$ to $e=5$ and $e=7$ to $e=11$, as well as the transition from $e=11$ at degree $k$ to $e=1$ at degree $k+1$. Equation \eqref{eq:common-sturm} gives a single Sturm-Liouville proof for the full modular cycle and also treats the mixed transition $e=5$ to $e=7$ directly.

\begin{theorem}\label{thm:complete-interlace}
For any two consecutive reducible parameters for which both non-elliptic pole divisors are nonempty, the two divisors strictly interlace on $\mathcal A$. More precisely, for every $k\ge1$,
\begin{align}
0&<x_{k,1}^{(1)}<x_{k,1}^{(5)}<x_{k,2}^{(1)}
<\cdots<x_{k,k}^{(1)}<x_{k,k}^{(5)}<1,\label{eq:int15}\\
0&<x_{k,1}^{(7)}<x_{k,1}^{(5)}<x_{k,2}^{(7)}
<\cdots<x_{k,k}^{(7)}<x_{k,k}^{(5)}<1,\label{eq:int57}\\
0&<x_{k,1}^{(7)}<x_{k,1}^{(11)}<x_{k,2}^{(7)}
<\cdots<x_{k,k}^{(7)}<x_{k,k}^{(11)}<1,\label{eq:int711}\\
0&<x_{k+1,1}^{(1)}<x_{k,1}^{(11)}<x_{k+1,2}^{(1)}
<\cdots<x_{k,k}^{(11)}<x_{k+1,k+1}^{(1)}<1.\label{eq:int111}
\end{align}
The same inequalities hold for the corresponding representatives on $\mathcal A$, with the order induced by the orientation from $\rho$ to $i$.
\end{theorem}

\begin{proof}
By Proposition~\ref{prop:common-sturm}, all four transformed polynomials satisfy
\[
(pU')'+\Lambda(r)wU=0,
\qquad
\Lambda(r)=\frac{r^2}{4}-\frac1{144},
\]
and $\Lambda(r)$ is strictly increasing for $r>0$. The consecutive parameters are
\[
2k+\frac16,\qquad
2k+\frac56,\qquad
2k+\frac76,\qquad
2k+\frac{11}{6},\qquad
2(k+1)+\frac16.
\]
On compact subintervals of $(0,1)$ the usual Sturm comparison theorem applies. Thus, when the spectral parameter is increased, there is a zero of the upper solution between any two consecutive interior zeros of the lower solution. It remains to control the intervals adjacent to $0$ and $1$.

For the transition $e=1$ to $e=5$, take
\[
U=Q_k^{-1/3,-1/2},\qquad
V=x^{1/3}Q_k^{1/3,-1/2}.
\]
Suppose that there is no zero of $V$ after the last zero $a$ of $U$, and choose signs so that $U,V>0$ on $(a,1)$. Then $U'(a)>0$, so
\[
W[U,V](a)=-p(a)V(a)U'(a)<0.
\]
Since $\Lambda_V>\Lambda_U$, \eqref{eq:wronskian-derivative} shows that $W[U,V]$ is strictly decreasing. Both solutions have nonzero finite limits at $1$, whereas Lemma~\ref{lem:endpoint-W} gives $W[U,V](1^-)=0$. This is impossible. Counting the $k$ zeros gives \eqref{eq:int15}.

For the transition $e=5$ to $e=7$, let
\[
U=x^{1/3}Q_k^{1/3,-1/2},\qquad
V=(1-x)^{1/2}Q_k^{-1/3,1/2}.
\]
Assume that $V$ has no zero before the first interior zero $a$ of $U$, and choose $U,V>0$ on $(0,a)$. Near $0$,
\[
U(x)=A x^{1/3}(1+O(x)),\qquad V(x)=B(1+O(x)),
\]
with $A,B>0$. Lemma~\ref{lem:endpoint-W} gives
\[
W[U,V](0^+)=-\frac{AB}{3}<0.
\]
Since $U'(a)<0$,
\[
W[U,V](a)=-p(a)V(a)U'(a)>0,
\]
contradicting the strict decrease of the Wronskian. Sturm comparison supplies the remaining zeros, and \eqref{eq:int57} follows.

For the transition $e=7$ to $e=11$, comparison gives a zero of the upper solution between each pair of consecutive zeros of the lower solution. If there were no additional upper zero after the last lower zero $a$, choose both solutions positive on $(a,1)$. Then $W(a)<0$. Both solutions vanish as a positive constant times $(1-x)^{1/2}$ at $1$, so Lemma~\ref{lem:endpoint-W} gives $W(1^-)=0$, again contradicting strict decrease. This proves \eqref{eq:int711}.

For the transition $e=11$ to the next $e=1$, take
\[
U=U_{k,2,2},\qquad V=U_{k+1,0,0}.
\]
Near $0$, the same argument as in the $e=5$ to $e=7$ case gives $W[U,V](0^+)<0$, while $W$ is positive at the first zero of $U$ if $V$ has no earlier zero. Hence $V$ has a zero in the first end interval. Near $1$, suppose that $V$ has no zero after the last zero $b$ of $U$, and choose $U,V>0$ on $(b,1)$. Then $U'(b)>0$, so $W(b)<0$. Here
\[
U(x)=A(1-x)^{1/2}(1+O(1-x)),\qquad
V(x)=B(1+O(1-x)),
\]
with $A,B>0$, and Lemma~\ref{lem:endpoint-W} gives
\[
W[U,V](1^-)=\frac{AB}{2}>0,
\]
which is impossible because $W$ is strictly decreasing. Sturm comparison gives one zero in every remaining interval. This proves \eqref{eq:int111}. The inequalities on $\mathcal A$ follow from the monotonicity of $J$.
\end{proof}

\section{Explicit solutions and an arithmetic application}

For fixed $k,\varepsilon,\delta$, let $\alpha,\beta$ be as in \eqref{eq:alphabeta} and put
\begin{equation}\label{eq:Fcanonical}
F_{k,\varepsilon,\delta}(z)=
\frac{\eta(z)^{4+8\varepsilon+12\delta}}
{E_4(z)^\varepsilon E_6(z)^\delta
\,Q_k^{\alpha,\beta}(J(z))^2}.
\end{equation}
By Corollary~\ref{cor:residue-from-ode}, the meromorphic differential $F_{k,\varepsilon,\delta}(z)\dd z$ has zero residue at every pole, including the elliptic poles when they occur. It therefore admits a global meromorphic primitive on $\Hh$.

Choose a point $z_0\in\Hh$ away from the poles and define
\begin{equation}\label{eq:hcanonical}
h_{k,\varepsilon,\delta}(z)
=\int_{z_0}^{z}F_{k,\varepsilon,\delta}(w)\dd w,
\end{equation}
where the primitive is understood by analytic continuation across paths avoiding the poles. Zero residues make it single-valued and meromorphic.

\begin{theorem}\label{thm:explicit}
Let
\[
r=2k+\frac{1+2\varepsilon+3\delta}{6}.
\]
Then $h_{k,\varepsilon,\delta}$ satisfies
\[
\{h,z\}=2\pi^2r^2E_4(z).
\]
A fundamental system for
\[
y''+\pi^2r^2E_4(z)y=0
\]
is
\begin{align}
y_1(z)
&=\frac{E_4(z)^{\varepsilon/2}E_6(z)^{\delta/2}
Q_k^{\alpha,\beta}(J(z))}
{\eta(z)^{2+4\varepsilon+6\delta}},\label{eq:y1}\\
y_2(z)&=h_{k,\varepsilon,\delta}(z)y_1(z).\label{eq:y2}
\end{align}
Up to post-composition of $h$ by a M\"obius transformation, this gives the solution of the Schwarzian equation for the corresponding reducible parameter.
\end{theorem}

\begin{proof}
Put $\ell=6r$. For $\ell\ge5$, Theorem~\ref{thm:KZ-pole-poly} and Proposition~\ref{prop:eta-normalization} show that the function $y_1$ in \eqref{eq:y1} satisfies
\[
y_1''+\pi^2r^2E_4y_1=0.
\]
The exceptional value $\ell=1$ corresponds to $y_1=\eta^{-2}$. In that case $D\log\eta=E_2/24$ and $DE_2=(E_2^2-E_4)/12$ give the same equation directly.

From \eqref{eq:Fcanonical},
\[
F_{k,\varepsilon,\delta}=y_1^{-2},\qquad h'=y_1^{-2}.
\]
Set $y_2=hy_1$. Then
\[
y_2''+\pi^2r^2E_4y_2=h''y_1+2h'y_1',
\]
and the right-hand side vanishes because $h'=y_1^{-2}$ implies $h''=-2y_1^{-3}y_1'$. Moreover,
\[
W(y_1,y_2)=y_1y_2'-y_1'y_2=y_1^2h'=1,
\]
so $y_1$ and $y_2$ are independent. Since $h=y_2/y_1$, Lemma~\ref{lem:schwarzian} gives
\[
\{h,z\}=2\pi^2r^2E_4(z).
\]
Any other quotient of a fundamental system differs from $h$ by a M\"obius transformation.
\end{proof}

The $J$-values of the pole representatives are algebraic, but the representatives on the elliptic arc are not.

\begin{theorem}\label{thm:transcendence}
For every $k\ge1$, every $(\varepsilon,\delta)\in\{0,2\}^2$, and every $1\le j\le k$, the point $w_{k,j}$ defined by \eqref{eq:wdef} is transcendental.
\end{theorem}

\begin{proof}
Put $\ell=12k+1+2\varepsilon+3\delta$. By Proposition~\ref{prop:pole-poly},
\[
j(w_{k,j})=1728x_{k,j}
\]
is an algebraic root of $\mathcal P_\ell$, and every root of $\mathcal P_\ell$ lies in $(0,1728)$. Suppose that $w_{k,j}$ were algebraic. Schneider's theorem then implies that $w_{k,j}$ is imaginary quadratic; see \cite{schneider}. Hence $j(w_{k,j})$ is a singular modulus.

Let $D<0$ be the discriminant of its imaginary quadratic order. The ring class polynomial $H_D(X)$ is the minimal polynomial over $\Q$ of a singular modulus of discriminant $D$; its roots are the $j$-values attached to the proper ideal classes of that order, see \cite{cox}. When the order is maximal this is the usual Hilbert class polynomial. Since the minimal polynomial of $j(w_{k,j})$ divides $\mathcal P_\ell$, all roots of $H_D$ would have to lie in $(0,1728)$.

The principal ideal class gives a conjugate outside this interval. If $D\equiv0\pmod4$, it is represented by
\[
\tau_D=\frac{i\sqrt{|D|}}{2}.
\]
For $D<-4$, one has $\operatorname{Im}\tau_D>1$. Along the imaginary boundary of the standard fundamental domain, $j(iy)$ is real and satisfies $j(iy)>1728$ for $y>1$. Hence $j(\tau_D)>1728$.

If $D\equiv1\pmod4$, the principal class is represented by
\[
\tau_D=\frac{-1+i\sqrt{|D|}}{2}.
\]
For $D<-3$, this point lies on the vertical boundary above $\rho$. On that boundary $j(-1/2+iy)$ is real and negative for $y>\sqrt3/2$, so $j(\tau_D)<0$.

The remaining discriminants $D=-3$ and $D=-4$ give $j=0$ and $j=1728$, respectively. Both are excluded because the roots of $\mathcal P_\ell$ lie strictly between $0$ and $1728$. This contradiction proves that $w_{k,j}$ is transcendental.
\end{proof}

\begin{corollary}\label{cor:noncm}
The elliptic curve represented by $w_{k,j}$ has algebraic $j$-invariant and no complex multiplication. The Jacobi pole configuration therefore gives algebraic points of the coarse modular curve whose representatives on $\mathcal A$ are transcendental.
\end{corollary}

\begin{corollary}\label{cor:ss-reduction-ec}
Assume that $\ell=12k+1+2\varepsilon+3\delta$ is prime. Let $K$ be a number field containing one of the pole invariants
\[
j_{k,j}=1728x_{k,j}.
\]
At every place of $K$ above $\ell$, the number $j_{k,j}$ is integral. After a finite extension of the corresponding local field, an elliptic curve with $j$-invariant $j_{k,j}$ acquires good supersingular reduction.
\end{corollary}

\begin{proof}
Theorem~\ref{thm:ss-lift} shows that $\mathcal P_\ell$ is monic with coefficients integral at $\ell$. Hence every root $j_{k,j}$ is integral at every place of its field of definition above $\ell$, and its reduction is a root of the non-elliptic supersingular polynomial $s_\ell$.

Let $K_v$ be the corresponding local field. The integrality of $j_{k,j}$ is the criterion for potential good reduction of an elliptic curve with this $j$-invariant; see \cite{silverman}. After a finite extension $L/K_v$, choose a model with good reduction. The $j$-invariant of the special fiber is the reduction of $j_{k,j}$, since $j$ is unchanged by base extension and specializes for a good model. By \eqref{eq:pole-ss-cong}, this reduced invariant is a root of $s_\ell$, so the special fiber is supersingular.
\end{proof}

\section{Asymptotic geometry of the pole configurations}

We keep $(\varepsilon,\delta)$ fixed and let $k$ tend to infinity. Write
\[
0<x_{k,1}<\cdots<x_{k,k}<1
\]
for the Jacobi zeros and let $w_{k,j}\in\mathcal A$ be their modular representatives.

For fixed Jacobi parameters, the zero counting measures converge to the arcsine distribution. Its pushforward under the inverse of $J$ gives the limiting measure on the modular arc.

\begin{theorem}\label{thm:equidistribution}
Let $\tau:[0,1]\to\mathcal A$ be the inverse of the restriction of $J$ to the oriented elliptic arc. For every continuous function $\Phi$ on $\mathcal A$,
\begin{equation}\label{eq:limitmeasure}
\lim_{k\to\infty}\frac1k\sum_{j=1}^k\Phi(w_{k,j})
=\frac1\pi\int_0^1
\frac{\Phi(\tau(x))}{\sqrt{x(1-x)}}\dd x.
\end{equation}
In particular, the pole representatives become dense in $\mathcal A$.
\end{theorem}

\begin{proof}
For fixed $\alpha,\beta>-1$, the normalized zero counting measures of the Jacobi polynomials converge weakly to the arcsine measure on $[0,1]$; see \cite{sze}. Since $w_{k,j}=\tau(x_{k,j})$ and $\tau$ is continuous on $[0,1]$, \eqref{eq:limitmeasure} follows by pushforward of measures.
\end{proof}

The endpoint scales reflect the ramification of $J$.

\begin{theorem}\label{thm:edge}
Let $w_{k,1}$ be the pole nearest to $\rho$ and $w_{k,k}$ the pole nearest to $i$. Then
\begin{equation}\label{eq:edge}
|w_{k,1}-\rho|\asymp k^{-2/3},
\qquad
|w_{k,k}-i|\asymp k^{-1}.
\end{equation}
\end{theorem}

\begin{proof}
For fixed Jacobi parameters $\alpha,\beta>-1$, the extreme zeros satisfy
\[
x_{k,1}\asymp k^{-2},
\qquad
1-x_{k,k}\asymp k^{-2};
\]
more precise Bessel-type asymptotic expansions are classical, see \cite{sze}.

At $\rho$, the Eisenstein series $E_4$ has a simple zero while $E_6(\rho)\ne0$. From \eqref{eq:J},
\[
J(z)=c_\rho(z-\rho)^3+O((z-\rho)^4),
\qquad c_\rho\ne0.
\]
Hence
\[
|w_{k,1}-\rho|\asymp x_{k,1}^{1/3}\asymp k^{-2/3}.
\]
At $i$, the form $E_6$ has a simple zero and $E_4(i)\ne0$, so
\[
1-J(z)=c_i(z-i)^2+O((z-i)^3),
\qquad c_i\ne0.
\]
Therefore
\[
|w_{k,k}-i|\asymp (1-x_{k,k})^{1/2}\asymp k^{-1}.
\]
\end{proof}

Both extreme Jacobi zeros are at distance of order $k^{-2}$ from the endpoints of $[0,1]$. The exponents $2/3$ and $1$ in \eqref{eq:edge} arise after pullback to the modular curve, where $J$ has ramification orders $3$ at $\rho$ and $2$ at $i$.

\section{Low-degree configurations}

The first primitive integer multiples of the monic polynomials $Q_k^{\alpha,\beta}$ are listed below, ordered by the residue class
\[
e=1+2\varepsilon+3\delta\pmod{12}.
\]

\begin{center}
\begin{tabular}{c c c c c}
\toprule
$(\varepsilon,\delta)$ & $e$ & $k=1$ & $k=2$ & $k=3$ \\
\midrule
$(0,0)$ & $1$
& $7x-4$
& $247x^2-260x+40$
& $2945x^3-4560x^2+1824x-128$ \\
$(2,0)$ & $5$
& $11x-8$
& $391x^2-476x+112$
& $4669x^3-8004x^2+3864x-448$ \\
$(0,2)$ & $7$
& $13x-4$
& $95x^2-76x+8$
& $5735x^3-7440x^2+2400x-128$ \\
$(2,2)$ & $11$
& $17x-8$
& $667x^2-644x+112$
& $1189x^3-1740x^2+696x-64$ \\
\bottomrule
\end{tabular}
\end{center}

For $k=1$ the residue system is linear and its unique solution is
\[
x_{1,1}=\frac{\alpha+1}{\alpha+\beta+2}.
\]
If $\ell=12k+e$ is prime, Theorem~\ref{thm:ss-lift} shows that, after $X=1728x$, the corresponding polynomial reduces modulo $\ell$ to the non-elliptic supersingular polynomial.


\begin{thebibliography}{99}

\bibitem{apostol}
T.~M. Apostol,
\emph{Modular Functions and Dirichlet Series in Number Theory},
2nd ed., Graduate Texts in Mathematics 41, Springer, New York, 1990.

\bibitem{babagranath}
S.~Baba and H.~Granath,
Orthogonal systems of modular forms and supersingular polynomials,
\emph{Int. J. Number Theory} \textbf{7} (2011), 249--259.

\bibitem{ram2}
K.~Besrour and A.~Sebbar,
Hypergeometric solutions to Schwarzian equations,
\emph{Ramanujan J.} \textbf{65} (2024), 1113--1121.

\bibitem{brillhartmorton}
J.~Brillhart and P.~Morton,
Class numbers of quadratic fields, Hasse invariants of elliptic curves, and the supersingular polynomial,
\emph{J. Number Theory} \textbf{106} (2004), 79--111.

\bibitem{cox}
D.~A. Cox,
\emph{Primes of the Form $x^2+ny^2$},
2nd ed., Pure and Applied Mathematics, Wiley, Hoboken, 2013.

\bibitem{cullinanhajir}
J.~Cullinan and F.~Hajir,
Algebraic properties of Kaneko-Zagier lifts of supersingular polynomials,
\emph{Proc. Amer. Math. Soc.} \textbf{145} (2017), 2291--2304.

\bibitem{driver}
K.~Driver, K.~Jordaan and N.~Mbuyi,
Interlacing of the zeros of Jacobi polynomials with different parameters,
\emph{Numer. Algorithms} \textbf{49} (2008), 143--152.

\bibitem{duncano}
J.~F.~R. Duncan and K.~Ono,
The Jack Daniels problem,
\emph{J. Number Theory} \textbf{161} (2016), 230--239.

\bibitem{forrester}
P.~J. Forrester,
\emph{Log-Gases and Random Matrices},
London Mathematical Society Monographs Series 34,
Princeton University Press, Princeton, 2010.

\bibitem{kz}
M.~Kaneko and D.~Zagier,
Supersingular $j$-invariants, hypergeometric series, and Atkin's orthogonal polynomials,
in \emph{Computational Perspectives on Number Theory},
AMS/IP Stud. Adv. Math. 7, American Mathematical Society, Providence, 1998, 97--126.

\bibitem{morton}
P.~Morton,
The Hasse invariant of the Tate normal form $E_5$ and the class number of $\Q(\sqrt{-5\ell})$,
\emph{J. Number Theory} \textbf{227} (2021), 94--143.

\bibitem{nakaya}
T.~Nakaya,
On modular solutions of certain modular differential equation and supersingular polynomials,
\emph{Ramanujan J.} \textbf{48} (2019), 13--20.

\bibitem{ogg}
A.~P. Ogg,
Automorphismes de courbes modulaires,
\emph{S\'eminaire Delange--Pisot--Poitou, Th\'eorie des nombres}
\textbf{16} (1974--1975), Exp.~No.~7, 8 pp.

\bibitem{forum}
H.~Saber and A.~Sebbar,
Automorphic Schwarzian equations,
\emph{Forum Math.} \textbf{32} (2020), 1621--1636.

\bibitem{ram1}
H.~Saber and A.~Sebbar,
Automorphic Schwarzian equations and integrals of weight $2$ forms,
\emph{Ramanujan J.} \textbf{57} (2022), 551--568.

\bibitem{jmaa}
H.~Saber and A.~Sebbar,
Equivariant solutions to modular Schwarzian equations,
\emph{J. Math. Anal. Appl.} \textbf{508} (2022), Paper No.~125887.

\bibitem{revista}
H.~Saber and A.~Sebbar,
Modular differential equations and algebraic systems,
\emph{Rev. R. Acad. Cienc. Exactas F\'is. Nat. Ser. A Mat. RACSAM}
\textbf{118} (2024), Paper No.~100.

\bibitem{schneider}
T.~Schneider,
\emph{Einf\"uhrung in die transzendenten Zahlen},
Springer, Berlin, 1957.

\bibitem{silverman}
J.~H. Silverman,
\emph{Advanced Topics in the Arithmetic of Elliptic Curves},
Graduate Texts in Mathematics 151, Springer, New York, 1994.

\bibitem{steinerberger}
S.~Steinerberger,
Electrostatic interpretation of zeros of orthogonal polynomials,
\emph{Proc. Amer. Math. Soc.} \textbf{146} (2018), 5323--5331.

\bibitem{sze}
G.~Szeg\H{o},
\emph{Orthogonal Polynomials},
4th ed., American Mathematical Society Colloquium Publications 23,
American Mathematical Society, Providence, 1975.

\bibitem{ittersum}
J.-W.~van Ittersum, G.~Oberdieck and A.~Pixton,
Gromov-Witten theory of K3 surfaces and a Kaneko-Zagier equation for Jacobi forms,
\emph{Selecta Math. (N.S.)} \textbf{27} (2021), Paper No.~64.

\end{thebibliography}
\end{document}